# Historical Changes in the Concepts of Number, Mathematics and Number Theory

Nicola Graves-Gregory
nicolasarasvati@gmail.com

*This essay traces the history of three interconnected strands: changes in the concept of number; in the nature and importance of arithmetike (αριθμητικη), the study of the qualities of number, which evolved into number theory; and in the nature of mathematics itself, from early Greek mathematics to the 20th century.*

*These were embedded in philosophical shifts, from the classical Greek ontologies through increasing pragmatism to formalism and logical positivism. Given Gödel's demonstration of the limitations of the latter as a foundation for mathematics, this essay explores phenomenology and Lakatosian ideas, which together offer a more sound epistemological and ontological basis for mathematics and a methodology for mathematical development.*

*The question also then arises of the possible resurrection of earlier, neglected mathematical projects, including widening the domain of number theory to include integer qualities revealed in the growth of mathematics in general, which has predominantly been the growth of quantitative mathematics or logistike (λογιστικη), the complement of arithmetike in classical Greece.*

## 1. INTRODUCTION

There were originally two main motivations in this research. One was a concern with the phenomenon of the integers. The intention was to examine how and why attempts to investigate the qualities of the integers gradually lost importance after the discovery of the existence of irrationals, and to examine what meaning the essential concern of the Pythagoreans with the integers might have for us today. This led on to considerations of the changes in the concept of number.

The other motivation was a dissatisfaction with the predominant trends in current history and philosophy of mathematics, in that they fail to adopt a critical perspective on present-day mathematics and are thus unable to play an active part in determining the telos of mathematical development, i.e. to offer constructive suggestions as to possible and desirable directions for mathematics. A form of synthesis: historical philosophy or philosophical history (as suggested by Lakatos and exemplified in his *Proofs and Refutations*[1]) might offer a real possibility of such a perspective. Certainly contemporary mathematical histories are not philosophically oriented in a presently active sense, and contemporary mathematical philosophies define themselves ahistorically[2]

This investigation, whilst beginning to examine these two initially distinct topics, hopefully indicates that they are, in fact, linked: that the attitudes and implicit assumptions which underlie current historical and philosophical perspectives are rooted in those that have governed the development of our conceptualisation of number and the perception of the integers.

Before going into more detail in justifying and amplifying the above statements and examining their implications, here is the outline of an argument for a historical philosophical perspective.

---

1 Lakatos (1964)
2 Since this writing Philip Kitcher (1983) has begun to approach this area.

1. Mathematical concepts and constructs are historically situated[3].

This statement does not conflict with a Platonist view of mathematics; it concerns mathematical means, not its telos.

To deny this statement would imply the belief that there had been no real change in mathematics (only a change in language), not a view generally maintained by mathematicians or metamathematicians.

2. We have choice as to how mathematics is to develop.

The contrary belief implies some kind of determinism.

Either a) behaviourist, i.e. we are out of our own control.

or b) mystical platonist, i.e. mathematical truths are absolute and eternal; mathematicians do not choose what they do mathematically, they can only try to follow their intuition as a guide to the truth[4].

If we accept 1 and 2, i.e. do not accept 2a) or 2b), then we are presented with the question:

3. In what directions do we think it desirable and possible for mathematics to develop?

This contains two questions:

3a) What are our criteria for deciding in what directions we wish mathematics to develop?

3b) What methods, modes, alternative aims, are possible for mathematics?

This problematic is ignored in contemporary mathematical philosophy, in that none of the mathematical schools concern themselves directly with the development of mathematics.

Since this problematic has not been explicitly, directly confronted before, we are not clear about the directions and aims implied in current mathematical praxis. So we have a preliminary question.

4) How can we clarify the intentional[5] meaning structure contained in contemporary mathematics?

Question 3a) obviously involves many questions including the social role of mathematics, one which has mostly been avoided in the present day although it was considered in the Greek era[6]. A full examination of these issues is outside the remit of this paper.

---

[3] This has been abundantly demonstrated by Becker (1927), Lakatos (1963/4), Fisher (1966/7), Dessanti (1968), Foucault (1970), Grabiner (1974) amongst others.
[4] See Gödel (1944).
[5] This is a phenomenological term cf. Husserl (1970).
[6] See example Plato's Republic (Part 8, Book VII) and Farrington (1946) for his reference to the 8th book of Plutarch's Dinner Table Discussions where Lycurgus is mentioned as allowing geometry but not arithmetic as a study in Sparta, because of the preferable political implications of the former.

Questions 3b) and 4) present an extremely difficult task since, as Husserl points out, our language, understanding and perception are permeated with the 'sedimentation' of history[7]. With any acquisition of knowledge, there is a forgetting of what our world was like before that. We all <u>learned</u> to see physical objects. What was seeing like before that learning?  We have learned to talk and so to think mainly in words. How did we think before we acquired language? This is a pertinent question for our present inquiry.

To attempt to extricate ourselves from this mesh, we must re-examine the history which created it: both the steps which led to the present mathematical configuration and the alternatives which were rejected. Although it is not possible to know definitively what was entailed in past possibilities, it is sufficient if we can begin

- A) to appreciate the nature of the contingent element in alternative mathematical projects,
- B) to understand why, in a given historical context, alternative projects were not adopted[8], and hence
- C) to understand the implications of the projects that were favoured.

Any significant progress here might enable us to tackle a perhaps more controversial project,

- D) to re-evaluate whether some previously rejected projects might have meaningful content and value for mathematics now.

We can find examples of such projects and approaches in the Greek and Renaissance/Reformation periods, which have been crucial in determining the content and modes of our present mathematics. In particular, the decline of Pythagorean number theory, in status and in living content, can be seen as one of the earliest examples of choice in the history of mathematics, and extremely valuable to examine.

Firstly let us look at the current state of the history of mathematics, in support of the contention that its lack of critical perspective produces serious weaknesses. It will be necessary to alternate between philosophical-historical concerns and those of number for a while. At times it may seem as though the windings and turnings of the argument have lost sight of their goal. But in fact the main concern throughout is the contemporary impasse of mathematical history and philosophy: and a concern to show that this impasse, when traced back to its historical origins, is deeply related to our relationship with the integers.

## 2. MATHEMATICAL HISTORY

Most research into the history of mathematics has attempted to understand the past from the perspective of contemporary formulations and concepts. This approach is valuable in picking out threads of continuity; but it also leads to misinterpretations and distortions of the past. A misleading emphasis is given to conceptual developments which can be seen as steps leading to present formulations, whilst those perspectives which in fact do not flow into the present are treated, in effect, either as pre-mathematical or non-mathematical, as 'anticipations' or 'false trails', thus obscuring not only the internal

---

7 See Husserl (1946) and Klein (1940). Also Barfield (1949) in his essay, 'The Force of Habit' (particularly pp. 69 – 79) describes vividly the imprisonment we experience as a result of our culture-based mental habits. He also (like Pythagoras, Plato, Descartes and others) emphasises the need for self-transformation to realise the extent of our imprisonment and to begin to release ourselves.

8 See especially Fisher (1966/7)

phenomenology of the past, but also the real nature of mathematical development[9].

Such a Whig theory of mathematical history - one which regards present concepts as the logically inevitable apex of mathematical investigation - is basically ahistorical: it makes it hard or impossible to see current mathematics, like all its earlier equivalents, as a moment in an <u>ongoing</u> historical process where, by no means inevitable choices are constantly being made between alternative projects and paradigms[10]. It is for contingent, not absolute, reasons that various historical projects have gone to the wall; and there are certainly senses in which the array of problems tackled by mathematics is historically arbitrary[11].

Because the underlying, historiographical premises are not generally explicated, they are not readily seen. When one does articulate the approach which most historians of mathematics have adopted, it becomes apparent that, in fact, an aspect of deductive logic (an element of contemporary mathematical theory) is imposed upon history on two distinct levels. On one level, clearly, in the Whig perspective outlined, which interprets earlier views and concepts as 'goodies' or 'baddies', true or false, right or wrong, according to whether they are perceived as being close to, or far from current views and concepts. On another level (where the implicit extension of an intra-mathematical attitude is perhaps not so immediately obvious), in the perspective of historians (sometimes consciously attempting to avoid the Whig distortions) who seek the <u>true</u> history, the <u>real</u> reasons for events etc., and fail to come to terms with the historical situatedness of their own perspective.

Certainly, in both these cases the picture is shaded, not sharp black-and-white, a multivalent truth function rather than a bivalent one. But the fundamental point is that a measure (albeit fuzzy) is imposed upon history: in the first case it is actively imposed upon historical phenomena; in the second, it plays a more passive role with respect to historical interpretations.

The historical roots of the black/white view could be said to lie in Parmenides' philosophy. If we do not accept the Parmenidean conclusion, but rather believe that change is real, then we see that,
As D.L.Miller says, 'the emergent gives rise to a new perspective, a new past'[12]; what seems plausible or important in a historical explanation changes as our mental frameworks change. This does not mean that history must necessarily act as a passive support for the status quo; it can provide a source of alternative perspectives, revealing hidden assumptions in those currently held. On examining the paths that led to our conceptualisations we may find places where the ideas that prevailed, did so not by transcending the former contradictions, but for contingent reasons, and re-examination of these issues may lead to new resolutions, other possible directions.

The 'sedimentation' of history permeates our attitudes and understandings to such an extent that reality is now seen as having an exact mathematical nature. That this was a superstructure imposed as a hypothesis, has been forgotten because the hypothesis has proven so successful in its intended sphere[13].

---

9 Bourbaki (1969) is one of the clearest examples of this Whig approach to the history of mathematics. See Unguru (1975) for a description of the results of such a perspective on the history of mathematics.
10 See Kuhn (1962), particularly c.11.
11 Wilder (1974) attempts to come to terms with this issue. See also Unguru (1975).
12 D.L.Miller (1948).
13 See Husserl (1970) and Barfield (1957). The basic fallacy underlying this (usually unconscious) assumption is pinpointed by a Sufi saying quoted by Keith Critchlow at a Wrekin Trust talk in June 1986, "Conclusive is not necessarily exclusive". The fact that one explanation of a phenomenon works well does not mean that it is the only possible explanation (that it is true

The technical usefulness of our current mode of mathematics cannot be denied. It has provided a wide variety of models that can accommodate quantitative change (one of the major advances from the Greek stage) and so serve for predictive scientific theories; but this was the original intention. Although it may be argued that present day pure mathematics does not share this goal explicitly, nevertheless it has developed within conceptualisations defined by scientific terms, and even, at a distance, by technological concerns: the demand for the calculus which arose from ballistics problems is just one example[14].

Obviously we have to examine more closely the areas and ways in which, for instance, the philosophical requirements of the Greeks are embedded in mathematical concepts and attitudes[15], as well as the metaphor transference (e.g. from the predominantly mechanical subject matter of the Renaissance) mediated through the seemingly neutral, abstract mathematical function.

What is key is to recognise that the success which quantitative mathematics and science have achieved in their chosen direction in no way validates its being the only direction possible for advance. We can draw an analogy with someone who wishes to leave a town: s/he is free to go in any direction; having journeyed, s/he may, at any time, measure her/his progress in terms of her/his distance from the starting point, but this offers no means for qualitative comparison between the place that s/he has reached and the other alternatives.

There is the possibility of approaching such a qualitative evaluation if we try to recover the mathematical problems of the past as they were formulated and understood in their own context, examining not only the ideas that have survived into the present, but also those that have been abandoned. We can then begin to understand the dynamics of mathematical development and obtain a more critical perspective on contemporary mathematics. By attempting to appreciate the meanings and implications of concepts and attitudes which persisted into the present mathematical corpus, we begin to have a context in which to discern the underlying intentional meaning structures as well as the possibility of examining the implications which alternative directions might have for us now.

This task of extricating ourselves, distancing ourselves from our ideological context is an orientation; the phenomenological '*epoché*' does not admit an absolute consummation, since that, like the historical perspectives already discussed, would imply transcending our historical situation. This is the goal that Klein and Husserl, in fact, set themselves. For Klein the final task arising from the attempt to reactivate the 'sedimented history' of the 'exact' nature is 'the rediscovery of the prescientific world and its true origins'[16]. The situation envisaged is, like that of the Cartesian doubt, not a real doubt, but a pretence; we cannot actually put ourselves into the pre-knowledge situation. The <u>attempt</u> to come closer to understanding the motivations and meanings (implicit and explicit) contained in earlier mathematical decisions must be seen in the context of an attempt to understand our own position.

Husserl was concerned to uncover the essential history of mathematics and thus the essential nature of mathematics. For him, to understand something is a positive act, a living moment; so one of the basic questions which he asks about mathematical development is: how could it be possible to relive all the moments of understanding that are necessary in a mathematical proof in order to progress to further knowledge? He sees knowledge as real, certain only in the moment of its realisation ('*Verwirklichung*').

---

and all others are false). A simple example (nonetheless valid) is that of the picture (often seen in textbooks on the psychology of perception) of the profiles of two black heads facing each other, separated (and created) by the outline of a white vase. It is not a question of either/or; both views are valid.

14 See Koyré (1957)
15 See Wittgenstein (1956), Cameron (1970) and Bremer (1973).
16 Klein (1940)

He recognises that it is impossible simultaneously to realise all the elements in a proof, (which would entail realisation of the complete reduction), so he creates (ad hoc) another category of thought to cover this case in mathematics, namely thought which has the potentiality of being reactivated[17]. The use of reactivatable elements in a proof would guarantee the soundness (reactivatability) of the whole[18], so that a vast edifice of mathematical knowledge could be built up from certain basic elements, provided that each step in the construction satisfied the criterion of reactivatability.

He does not consider the question which is prime for Lakatos (which has been put by mathematicians through the ages), namely, how does one arrive at mathematical concepts or discover (or create) new theorems in the first place, a stage which is necessarily prior to any attempt at proof. Focussing on mathematical development rather than certainty, Lakatos advocated a different attitude to the concept of proof from the present norm, and from Husserl's idealisation which is in fact structurally the same as the norm but projected onto a deeper meaning level.

Rather than seeing the attempt to prove a theorem as the attempt to establish it beyond doubt, in which case the appearance of counterexamples is regarded primarily as a failure of the theorem (or sub-lemmas), Lakatos was concerned not with an absolute result, but with the process involved in attempts to prove a theorem. By setting out the various steps of the argument, by looking for local and global counterexamples (i.e. counterexamples to the sublemmas, and counterexamples to the principal theorem), the hidden assumptions in the concepts used may be revealed, opening the way to realisation of more general concepts underlying those initially considered. Proof attempts and counterexamples thereby act as two poles in a continual process of improving conjectures, refining our mathematical ideas.

Interestingly, the present day term, 'proof' replaced the earlier 'demonstration', which is a closer translation of the original Greek 'δεικνυμι' (*deiknumi*) meaning, 'I show'. 'Proof' not only lays greater claim to certainty, it also shifts the power balance. When someone *shows* you an argument, then you might or might not be convinced by it; if someone *proves* it, they are claiming authority. Furthermore the term 'proof' assumes that there is an objective truth or an objective rationality, independent of any subjective consciousness assessing the rational process. The term 'δεικνυμι' seems closer to Lakatos' interpretation[19].

In fact Husserl's initial analysis (when we ignore his ad hoc creation of 'reactivatability'), is a valid approach to a phenomenology of mathematics which complements Lakatos' suggested methodology. It is because it is not possible to carry out a simultaneous reactivation of all the steps involved in a mathematical argument, that all proofs are temporary and non-absolute. They are necessarily partial, since the limits of the definitions of concepts cannot become clear until they are seen to be contained in deeper, more general concepts[20]. Because it is impossible to apprehend a proof in its totality back to the basic premises, attempts to reactivate certain elements of a proof, in their context or in another wider context, can reveal new ideas that were obscured by implicit assumptions in the original proof.

---

17 Husserl (1970). See also the discussion of nomological and eidetic disciplines in Husserl (1913).
18 Husserl (1929).
19 Any information on when, where, how and why the term 'demonstration' was replaced by 'proof' will be gratefully received.
20 See Lakatos (1963/4) and Polya (1954). Lakatos demonstrates this process of increasing inclusion. See also Whiteside (1960/2). Einstein's gravitational theory contains Newton's. Gödel proved this common sense idea in mathematics. It is interesting and unfortunate that although the more inclusive concepts and theories do not contradict earlier more partial ones, yet they are usually seen at first as threatening.

By following the development of a geometrical problem into topological concepts, Lakatos showed that the process he recommended as a methodology was that which actually took place over the course of history, but unconsciously (and in fits and starts) since, in keeping with the prevalent attitude to proof, counterexamples were regarded as 'monsters', causing mathematical reactions which he characterised as 'monster-barring', 'monster-adjusting' etc., rather than stimulating careful re-examination of the steps of the proof to discover precisely which steps or concepts were called into question in each case.

Lakatos' suggestion was that this process of improving conjectures should be adopted consciously as a methodology of mathematical development. In effect, he argues that mathematics already has a methodology latent within present forms; it is only the rigid, static interpretation of these forms which has prevented perception of the essential, complementary dynamic that they contain. It requires only a fluid rather than a static attitude to proof in order to appropriate and activate this methodology (to work with it rather than against it).

Lakatos' critique focusses on attitudes to the concept of proof, which determine reactions to the occurrence of counterexamples, thus affecting the development of any given problematic in mathematics. This paper, like Lakatos, is concerned with underlying attitudes which have influenced the development of different modes of mathematics.

The way in which most mathematical historians regard as deviants earlier mathematical projects which do not translate or conform to the current, conventional mathematical ideas, is similar to the view of counter examples as monsters; the insidious extrapolation of the mathematical-logical view of contradiction from static to temporal phenomena manifests itself in the inability to come to terms with real change.

Szabó suggests that the reductio ad absurdum form of proof which replaced earlier more illustrative forms, was adopted from Parmenidean logic[21]. Whether or not this was the case (his argumentation is convincing) the two forms are effectively the same. This shift in the mode of proof resulted in a move away from the intuitive reasoning which led to the conjecture initially[22]. This particularly affected the proofs of number theoretical theorems[23], which is surely a factor in the impoverishment and degeneration of arithmetike, which accompanied the subsequent growth of geometry. We shall see that this was also the beginning of the notion of mathematics as being more essentially concerned with quantity than qualities.

Obviously these questions require a much fuller investigation. Nevertheless we can begin to see that at the very inception of the mathematico-logical method, which was to prove so powerful in the development **of** mathematics, it brought about changes which profoundly affected the telos of mathematics. The implications of some of these changes were only articulated much later. It was still later before some of the inherent limitations began to be realised and the question as to the limits of validity (and usefulness) of the method recurs[24].

---

21 Szabo (1969).

22 See Szabó (1969). See also Whiteside (1960/2) for a discussion and illustrations of this shift as it affected mathematical proofs from the 17$^{th}$ century into the 18$^{th}$. The question of the meaning of proof and the validity of proof now recurs with extra bite given the computer proof of the 4-colour theorem and the collaboration proof of the classification theorem for finite groups.

23 Szabó (1969).

24 For example the technical consequences of some of Kronecker's arguments (e.g. on Aristotelian logic and on the centrality to mathematics of the natural numbers) were only realised fully in the work of Brouwer and the intuitionists. It is ironic that Aristotle (c.350 B.C.) was clear about the limits of validity of his logic; it is those who have used it subsequently who have implicitly assumed it to have global validity.

Husserl and Lakatos are important, because they open up debates in little touched areas of mathematical philosophy: Husserl's discussion relates not only to epistemology but also to ontological issues in mathematics. Lakatos is concerned with the methodology of mathematical development. The latter's study, using historical research to offer a constructive critique of current mathematics and mathematical philosophy, also serves as a relevant example of the results possible through combining historical investigation and philosophical analysis. I think it can be shown that Husserl and Lakatos, via different paths, in fact enter one space, a philosophical area peculiar to mathematics with its unique relation to the life-world.

Now, let us return to the mathematical focus of this inquiry: the concept of number.

## 3. NUMBER THEORY: NOTES TOWARDS A SELECTIVE HISTORY

This section will deal with topics which concern the integers; but not directly with the integers themselves. As necessary background, therefore, here is a short eulogy to natural number.

Who cannot be intrigued by the dual nature of the integers? In the first place they are our archetype of a discrete well-ordering. (It is this aspect which is taken as their defining characteristic in attempts to 'found' the integers on various set theories, i.e. to construct set theoretic models that functionally approximate to the natural numbers.) They appear as a monotonous repetition of a single relationship ad infinitum: $1<2<3<4$ ...... very useful for counting sheep, but scarcely exciting.

And yet as we look closer, as our range of numerical operations expands, we perceive more and more complex structures generated by this seemingly banal series: more varied relationships between the numbers, in terms of which we begin to appreciate the integers as individuals with different characteristics. As the range of structures that we are able to identify increases in diversity and complexity, where we once saw undifferentiated extension, we now see a finer, more subtle web of interlaced and distinct entities. We then see this web to have been latent in our primary, deceptively simple sequence. When an ever finer grain emerges at every increase in our power of definition, the subtlety of the number series seeming always one step ahead of our subtlety, how can we avoid the induction of the image of the natural numbers as the eternal source and limit of our pattern-making[25]? With the primes, perhaps, as perpetual jokers, continually escaping the webs we weave.

Well, yes, it seems natural to be in love with the natural numbers; and, as will emerge in this essay, the real number line and the complex plane can be seen as natural extensions of the integers: extensions whose developments and properties seem at times, like the paths in Alice's looking-glass garden, to lead away from their source, but always finally return to it.

When we look at contemporary mathematical work on the integers, although various mathematicians have commented on the integers in this vein, very little actual mathematical work corresponds to this attitude[26]. There are two main areas of contemporary mathematical work on the integers: foundational studies and number theory.

In the former, the approach adopted to the integers is that it is necessary to found them in some form of set theory[27]. This results in a model of the integers as a class, taking their most obvious attribute, the discrete well-ordering, as the defining characteristic, and (since the stress is on the homogeneous aspect of their nature) reveals nothing about the complex interrelations of the integers (i.e. about the integers as individuals), but then that is not its aim.

---

25 This is what the Löwenheim-Skolem theorem says.
26 Or even work in the philosophy of mathematics: one exception is Becker (1927). [N.B. 2013 note: the mathematical landscape has begun to change with the growing importance of cryptography and the discovery of the deep importance of the mysterious Riemann hypothesis.]

Thus foundational work has a more clearly articulated aim than contemporary number theory: to ensure secure foundations for the mathematical edifice. It is perhaps partly for this reason that it has enjoyed a higher status; even though the strong form of this aim that the system should be consistent and complete, was shattered years ago by Gödel's theorems. Certainly foundational, mathematico-logical concerns continue to dominate mathematical philosophy. But we shall see that the driving force behind this emphasis is an attitude which has haunted mathematics and its meta-disciplines, mathematical history and philosophy since Parmenidean bivalent logic was incorporated into mathematics as a guarantor of certainty, Wittgenstein drew attention to it in 'Remarks on the Foundations of Mathematics', where he writes

'My aim is to alter the attitude to contradiction and to consistency proofs. Not to shew that this proof shews something unimportant. How could that be SO?'[28]

This attitude, which could be called 'contradiction-phobia', is a fundamental block that we are in a far better position to overcome now than in the Greek era.

Needless to say, the dominant concern of this paper is closer to the other main approach to the integers, via number theory, which takes the integers as given and concerns itself with the integers as individuals, their characteristics and interrelations. For this reason we shall examine how and why, since the time of the Pythagoreans, number theory gradually declined in status and lost sight of its original telos, thus losing touch with its living content and losing wholeness and coherence. We shall also examine whether there might be a way in which the essential number theoretical concerns of the Pythagoreans could be meaningful today.

For the Pythagoreans, number theory or arithmetike, was the basis of their metaphysical science[29]. Their monadology was an attempt to discover the relationships of the universe, which they originally believed could be described totally in terms of integers and ratios of integers. In their mathematical work they both extended the range of possible operations with numbers (theory of ratios, magnitudes in geometry etc.) and simultaneously objectified the properties of the individual numbers which emerged[30]. That is to say, they were concerned both to develop more complex structures which could accurately describe the complicated phenomena in the world, and to take note of the role which the integers played in these developments.

By naming the qualities of the integers (i.e. properties which are the complementary result of an operational pattern, e.g. triangular, square etc.), it may prove possible to perceive relationships between these properties themselves (i.e. another mathematical structure) as, of course, has proved to be the case with the properties which the Pythagoreans abstracted. These relationships between the Pythagorean-derived properties are still the foundation of number theory today. It is possible gradually to build up a more complete picture of each individual integer in terms of its properties; in this case, in whatever context one is dealing with a particular integer, one's awareness of its various attributes may yield simultaneously a new understanding about the nature of the context (thus revealing new possibilities) and about the nature of the integer, for example another attribute, or a meta-relation of properties.

---

[27] This attitude which has been shared in one form or another by all the great figures in modern philosophy of mathematics from Frege (for whom *Wertläufe* take the place of sets) onwards, has been criticised by Wang who has pointed out that no substantial reasons have been given for supposing that 'numbers evaporate while sets are rocks' (Wang, 1974 p.238).See also the more general criticisms of 'foundationalism' in Smith (1976) and the references given there.

[28] Wittgenstein (1956)

[29] The distinction made by Klein (1968) between number theory and arithmetike has some validity and leads to certain interesting questions but they are not immediately relevant to the present discussion.

[30] For example, as Aristotle (1942, Δ.14.1020a35 – b8) tells us, 'linear', 'plane' etc

We do not know if the Pythagoreans <u>consciously</u> adopted this twofold approach. In fact, whereas the current mathematical approach emphasises operational structures, the Pythagorean emphasis seems originally to have been to discover qualities of the integers and structural developments served as tools to this end. But the resultant of their work was a balance between these two aspects of development (the dependence of the hypostatisation on operational developments is not matched by the inverse necessity) until the discovery of the existence of irrationals.

There is a tendency to depreciate the Pythagoreans' concern with the qualities or forms ('eide', ειδη) of the integers as number mysticism, which is questionable. There are two interwoven strands in their approach: the attempt to associate different numbers and ratios with human characteristics, moral attributes etc.; the other we could retrospectively (anachronistically) call defining equivalence classes of integers (largely from the forms which arose out of the figurate representations). The former does not invalidate the latter, as the attempts of the Merton scholastics to quantify such entities as love etc. do not invalidate their work on dynamics.

The question of the possible importance of a mystical or metaphysical perspective in furthering the more limited mathematical praxis will not concern us now. This question arose with regards to the relation between Newton's alchemical researches and his accepted mathematical work[31]. What is clear is that the usual connotation of 'mystical', namely 'opposed to reason', does not apply to the Pythagoreans. Their transcendent telos is rooted in a materialist perspective: they are described by Aristotle, together with the physiologists, as being of the opinion that being extends no further than sense perception[32]. Aristotle also stresses that it is Plato who first makes number separable from objects of sense, whereas for the Pythagoreans "the monads have magnitude"[33].

The discovery that there were magnitudes that could not be expressed as integers or ratios of integers meant that the Pythagorean monadology was no longer tenable as a global philosophy. The original concept of number was maintained. The term 'number' ('arithmos', αριθμος) was used only for the integers greater than one. 'One' (the 'monas', μονας) was the unit and unity, and was considered the principle or beginning of number and as such had a different status from the numbers which it generated. 'Two' was sometimes similarly excluded from the realm of number; the reason for this stems from the Pythagorean identification of 'one' and the odd numbers with 'limited', opposing 'two' (or the 'dyad') and the even numbers as 'unlimited'[34]; but this perception of two was not observed so strictly. Fractions were

---

31 See Dobbs (1975) who shows that Newton's alchemical work preceded and was concurrent with his scientific research. It was previously assumed that he became an alchemist only later in his life when he had completed his scientifically acceptable work. This assumption made it possible to dismiss the alchemical work as simply the product of his dotage.
32 Aristotle, (1942, A8, 990a 3f)
33 Aristotl2, (1941, M6, 1080b 19f)
34 The importance of unity, oneness, wholeness, is beginning to be recognised again with the growth of holistic approaches to knowledge. In biology clearly the whole is greater than the sum of its parts. The argument here is that we need to begin to investigate 2-ness, 3-ness etc in similar not considered to be numbers, since the 'one', the source of numbers, was essentially indivisible; only ratios of integers were allowed in arithmetike.

Since the irrationals are not generated by the 'one', nor do they reduce to ratios of integers, they had no place in this system (they did not even have a justified operational framework until Eudoxus adapted the theory of proportion to this purpose), and they were incorporated into mathematics as geometrical magnitudes (not numbers). There was now a rigid ontological distinction between the objects of study of arithmetike and geometry: 'number', arithmos, is discrete, a multitude of indivisible units; 'magnitude' ('megethos', $\mu\varepsilon\gamma\varepsilon\theta o\varsigma$) is continuous, an infinitely divisible spatial measure[35].

From this time on, the development of geometry began to outstrip that of arithmetike. Before this, the figurate representation of numbers and the theory of ratios, meant that geometry and arithmetike were closer and shared developments; now their objects of study were separate. It might be argued that, in any case, number theory had already exhausted the possibilities for deriving terminology from geometric forms. This is a separate question; here, what is relevant is that number was excluded from the study of geometry. In the Renaissance when the distinction was bypassed (with no theoretical underpinning), the content of number theory was already determined and the living relationship between the study of the integers and other mathematical areas was not renewed.

From this time also, there was a shift from illustrative demonstration to more strictly logical proofs - relying more on reductio ad absurdum. Szabó shows that for number theoretical theorems, such indirect proofs were sometimes substituted unnecessarily (in terms of rigour) and perversely (in terms of intuitive value)[36]. By the time of Euclid's compilation, the appended diagrams served no useful purpose, representing discrete numbers by line segments.

Also in this period the study of number itself was split into two disciplines: arithmetike and logistike. The dividing line was never unequivocally established, but one fairly common factor in the various versions given is that arithmetike deals with the 'eide', forms, kinds, species, of number; whereas logistike deals with the 'hyle'($\dot{\upsilon}\lambda\eta$) of numbers, the quantity, the material, the amount that they represent[37]. The verbal roots of their mathematical meanings are respectively, 'arithmein', to count, and 'logizmein', to calculate.

Disregarding for the time being the original underlying reasons for this distinction, it is clear that it was uncritically accepted by the neoplatonists as a rigid separation (Plato's suggested refinement of a further classification into theoretical and practical areas of each was not followed). The attempts at clarifying the distinction consisted in trying out different formulations for it and adjusting the classifications of the existing mathematical subject classification accordingly, rather than questioning its basic premises. So the study of the natures of numbers and that of numerical operations were seen as separate rather than as complementary and mutually stimulating.

It is relevant that Diophantine analysis, which was vitally important for the growth of modern algebra and which, according to the Platonist distinctions, should have appeared as logistike, or theoretical logistic (when Vieta takes it up, he adopts the latter term) in fact appeared in his Arithmetica**,** i.e. Diophantus disregarded this distinction. He also moved away from the mainstream in that he did not use the Euclidean proof form and he accepted fractions as numbers.

---

ways to those indicated by M.-L, von Franz (1974), using the wealth of modern mathematical knowledge, mainly logistike, to uncover the corresponding arithmetike.
35 These issues are discussed in detail by Becker (1927), pp.129-148 and pp.199-213.
36 Szabó (1958), pp.118-120.
37 See Klein (1968) for more detail.

So we have now extricated three important factors contributing to the decline in status and/or content of number theory in the Greek era:

i) the rigid ontological distinction between number and magnitude,

ii) the emphasis on logical proof (particularly bivalent logic) as opposed to illustrative demonstration,

iii) the hypostatisation of the arithmetike/logistike distinction.

Obviously these points require further examination, but now we continue our whistle-stop history of number theory. In the Middle Ages interest in the integers was primarily in terms of numbers as religious or moral symbols; this is certainly of interest in some respects but it is not directly relevant to our inquiry at present. From the Renaissance until the 19th century number theory basically consisted of a range of seemingly rather disparate (albeit stimulating and important) problems, such as the question of the distribution of primes, Fermat's last theorem etc., that appeared to have arisen almost accidentally in the course of its history. It had become an area of mathematics that lacked an inner sense of direction and wholeness, having been content to assume problems which involved terms which figured in the Greek number theory. Its only claim to wholeness was the tenuous continuity with the Greek discipline, which it maintained by preserving the superficial content of their concern. One might say that number theory had petrified.

In the 19th century it gained coherence when Gauss (who considered number theory to be the queen of mathematics) published his *Disquisitiones Arithmeticae* which extended the notion of integers to include complex integers, laying the ground for algebraic number theory. Analysis began to be used in number theoretical proofs, giving rise to analytic number theory. But the emphasis continues to be on other branches of mathematics acting as investigatory or proof machinery with respect to number theoretical problems whose roots go back to Pythagorean arithmetike. Although analytic and algebraic insights have extended and deepened the structural vocabulary of number theory, it has not appropriated the expanded field of operations involving number (the extensions of logistike) as a potential source for developing its basic descriptive terminology for the integers themselves.

The phenomenon of the occurrence of particular integers in diverse mathematical fields is not examined for the significance they may have in terms of a nature or characteristic of the integer involved; even though, particularly in algebraic geometry and topology, such phenomena are increasing and it is sometimes necessary to use classical number theory in such proofs, still the converse approach is not adopted. Obviously this is now a formidable task, but results seem to be converging in this direction.

Although it would be worthwhile to consider some perspectives on arithmetike as the matrix for logistike, for now we shall make a preliminary investigation of some of the questions raised by the history of the extension of the number concept.

## 4. CHANGES IN THE CONCEPT OF NUMBER

We shall firstly look at the way in which the concept of number was extended in the European rebirth of mathematics, at some of the underlying attitudes and their implications. We have already briefly considered the Greek conceptualisation: their ontological distinction number : magnitude corresponded to the antinomy discrete : continuous. Before beginning to investigate some of the philosophical questions raised in connection with the changes in the concept of number, we will first go through a very summary history of the developments after the Greek era.

The Romans were primarily interested in practical results rather than theoretical mathematics and the continuing usage of their number system in the 'dark' ages (making multiplication and division extremely lengthy tasks) meant that there was little theoretical mathematical work in this period. Its rebirth was an important element in the phenomenon of the Renaissance. Theory was stimulated by the introduction of Arabic and Greek texts (the latter via Arabic translations). The gradual adoption of the Hindu-Arabic number system which contained a sign for zero and was a consistent place system, facilitated numerical operations. This notational development was key to the explosion of mathematical, scientific and economic developments that followed.

It was in the sphere of commerce that this number system was first introduced: the main concern was correct, convenient operation with numbers, not a theoretical foundation. Since fractions arose in simple numerical calculations, they were popularly considered numbers, i.e. the Greek discrete:continuous, number :magnitude distinction was bypassed.

The study of the Arabic texts gave rise to much interest in algebra, primarily in solving polynomial equations: the terms 'surd', 'absurd' or 'cossike' numbers were used variously for rationals, irrationals and what we would now call the variable terms of an equation. There was an implicit assumption that polynomial equations were determinate, that there was a definite numerical solution waiting to be discovered (or, as we would describe it, that the existing number field was closed under the operations used).

The use of the term 'number' for these cases was disputed, in the first place by the neoplatonists. The emergence of negative and imaginary solutions caused further confusion. There was in fact more difficulty in accepting negatives than irrationals[38] since the negatives lacked the framework which the Eudoxan theory of proportions supplied for the irrationals. From the general questioning as to the criteria for acceptance of candidates for numberhood, the dominant perspective which emerged was a pragmatic one: operations with the new number-like entities continued because they were useful. Stevin's position was one of the most coherent; amongst others he championed the decimal notation for fractions and consequently advocated the radical notion that *'Number is not all discontinuous quantity'*[39].

This was the vague beginning of the idea of a number line, a different kind of parallel between geometry and arithmetic from that of the Greeks. The existence of a symbol for zero was of vital importance in this. Stevin and later Wallis argued its equivalence to the geometric point, as opposed to the Greek equivalence of the 'monas' to the point. They were also both concerned to be explicit about the corresponding new status of 'one', that it should be considered a number, since, according to Wallis, it answers the question, 'How many?'.

In the 17th century the question of the conceptualisation of number was generally secondary to problematics of the operational developments. As a result of the more pragmatic attitude (revealing the beginnings of a formalist attitude to mathematics) negatives and imaginaries were used like other numbers in calculations because they ultimately rendered correct results, even though when they emerged themselves as solutions, these were regarded as meaningless[40].

These operations and the corresponding attitude were not validated until the l9th century, when Hamilton elaborated the consistent algebra of complex numbers, and negatives and imaginaries were accorded a more 'intuitive', visual meaning in the Gauss-Wessel representation of the complex plane.

In the 18th century there had already been attempts to prove the fundamental theorem of algebra which implies that no new candidate for numberhood could emerge from polynomial equations. Leibniz in 1682 had coined the term 'transcendental' for numbers that 'transcend the power of algebraic methods'. It was in the 19th century that first Liouville (in 1844) demonstrated that a certain serial form would yield transcendental numbers, and later Hermite (in 1873) showed the transcendence of e, and Lindemann (in 1882) that of pi.

---

38 See Becker (1937) and Kline (1972).
39 See Klein (1968)
40 This stance is similar to the Copenhagen interpretation of quantum mechanics.

Both irrationals and transcendentals slip through the dense mesh of decimal fractional notation for the number line and at the end of the 19th century, Cantor developed a theory of transfinite numbers (cardinal numbers of infinities) which included a proof that the order of infinity of the continuum was higher than that of the whole numbers (and rationals). The development of the theory of these new numbers opened up controversies that harked back to the Eleatic paradoxes. Also, interestingly, a new monadic structure emerges in that the transfinite numbers are discrete: no-one has yet managed to construct a set whose cardinal number lies between that of the natural numbers and that of the real numbers, and higher transfinite numbers are generated by considering the set of subsets of a transfinite set[41].

Now to return to some of the philosophical questions raised. First, before considering any of the more particular questions, when examining Greek mathematics we cannot ignore the context of the original mathematical concern, which takes us into contradictions whose interplay has been a vital element in the growth of the mathematical organism.

For the Pythagorean-Platonic perspective (from which, rather than from the other Greek schools of thought, our mathematics mostly derives) mathematics was intended primarily as a metaphysical discipline; the aim was to understand material reality as a manifestation of divine truth, not to control nature. The guiding principle was the order ('taxis', ταξις) of the whole; the foundation of their metaphysics was a belief in the eternal unity, the 'one', and they sought, through an understanding of changeless mathematical relations, to reach the highest truths pertaining to the eternal reality which transcends material reality[42].

For the Pythagoreans, the transcendental, absolute truth was immanent in material, phenomenal reality; for Plato it lay behind or above worldly reality; there was a separation, an abstraction. For Pythagoras mathematics contained the truth; for Plato (although there is some ambiguity about his position) it seems that mathematics was a step towards an appreciation of the truth. For both, perception of the truth necessitated a self-transformation: it is not nature that hides, it is our vision that is skew[43].

Since the time of the Greeks, mathematics has been identified with truth, but the nature of the truth sought has changed radically from being a revelation of transcendent reality, to an undeniable, empirical statement or fact[44], or even a system of tautologies. Paradoxically, the one constant element in the idea of truth has been the quality of being absolute, timeless, unchanging! The underlying issues here touch the shared boundary, the no-man's-land between epistemology and ontology.

---

[41] Paul Cohen's exciting proof of the independence of the continuum hypothesis and axiom of choice from the axioms of ZF set theory, does not change this.

[42] See Becker (1927), Cornford (1939) and Klein (1968). The shift in attitude is clear if we look at the original meaning of 'theory'; as Russell says (1971, p.52), "this was originally an Orphic word which Cornford (in 'From Religion to Philosophy') interprets as 'passionate, sympathetic contemplation'". In this state, he says, "the spectator is identified with the suffering God, dies in his death, and rises again in his new birth". For Pythagoras, the 'passionate, sympathetic contemplation' was intellectual and issued forth mathematical knowledge.

[43] The term 'mathematics' was coined in Pythagoras' school: 'ta mathemata' (τα μαθηματα) meaning 'those things which have been learned'. Given the religious nature of the Pythagorean school, it is certainly arguable that the fundamental learning concerned the spiritual development of the disciples and what we now consider to be the mathematics deriving from that school, were symbols of that spiritual or psychological development. It is interesting that the original etymological meaning of 'mathematics' is so open. This may account for one problem which has dogged mathematical philosophy for some centuries (although it has not always been acknowledged), namely that mathematics has been defined by extension not by intension (to return to a neglected but useful Aristotelian distinction).

The Pythagorean-Platonic orientation was primarily ontological. The Platonic dialectic was a process to allow an ascent to a vision of the form of the good, a transcendent, subjective, absolute certainty, which was essentially incommunicable to others who did not take part in the ascent. Proclus' neoplatonist exposition of an 'ascent from more partial to more universal understandings' by which 'we climb up to the very science of 'being' in so far as it is 'being''[45] was extremely influential in the Renaissance; but the telos behind this ideal of a science transcending other sciences was profoundly different from the original Platonic goal. Plato's vision of a unified metaphysical science was split. The subjective, metaphysical orientation was preserved in the opera of the alchemists; but these, for various reasons, including the danger of persecution, largely remained privatised and that mode ultimately disappeared. Recently such works are being re-evaluated in terms of their effects on the new, emergent science, the rival mode which triumphed.

It is, of course, the Cartesian philosophy which most clearly indicates the reversals of the classical Greek approach, as well as the continuities. Plato's dialectic works through accepted statements to attain a higher level of absolute, subjective certainty; with Descartes' epistemological orientation, he reverses this order to found his rationalism on subjective certainty of a low level, then, taking mathematical reasoning, the 'geometric' method as his model, he believes it possible to proceed via clear-cut self-evidences, to accumulate a higher level of knowledge that is explicit, articulatable and still absolutely certain[46].

In his *Discourse on the Method of Right Reasoning*, Descartes says, "those long chains of reasoning, simple and easy as they are, of which geometricians make use in order to arrive at the most difficult demonstrations, had caused me to imagine that all those things which fall under the cognisance of man might very likely be mutually related in the same fashion and provided only that we abstain from receiving anything as true which is not so, and always retain the order which is necessary in order to deduce the one conclusion from the other, there can be nothing so remote that we cannot reach to it, nor so recondite that we cannot discover it."[47]

Remarkably no one appears to have commented on the fact that Descartes, in his enthusiasm for the 'geometric' method which is deductive (i.e. a method to prove a theorem once it exists); mistakes it for an inductive method (i.e. a means for discovering or generating new results).In fact he intended his method to be used only after his course of meditations had been followed. This part of his teaching appears to have been totally ignored both by his disciples in his lifetime and subsequently.

As Husserl[48] points out, there is also a fundamental inconsistency between Descartes' radical starting-point, the epoché, and the rationalist system that he develops (the former destined ultimately to undermine the latter): his bracketing is incomplete; belief in the Galilean, mathematical book of nature is not submitted. Before following any further the theme of the nature of mathematical truth and certainty, attempting to unravel the cross-threadings of the empirical and transcendental, objective:subjective, relative:absolute etc., let us return to consider the development of the number concept and its relation to the more general conceptualisation of mathematics as a whole.

---

44 See Lewis Mumford (1970). It is interesting that the etymological roots of 'fact' and 'fiction' are so similar: 'fact' derives from Latin 'facere', to do or to make; 'fiction' from Latin 'fingere', to form or to fashion.

45 See Klein (1968)

46 Descartes (1973)

47 Descartes (1973), p.92

48 Husserl (1970) part ll, c.18, p.79 ff.

For the Pythagorean and Platonic perspectives, the question of the ontological status of number and mathematical concepts is of prime importance. In the original Pythagorean conceptualisation this status is clear: all is number, where number is discrete, heterogeneous, substance and form, quality and quantity. The existence of the irrationals, which made this early vision of a complete, mathematical universe untenable, was encompassed in Plato's hierarchy of levels of reality. The irrationals can be seen as magnitudes and part of the illusory, changeable, material world; the integers belong to the real, metaphysical, changeless realm of ideals, where the form of the good presides. The discrete:continuous, number : magnitude contradiction was contained, if not clarified, in an ontological distinction.

The beginnings of the modern conception of number can be seen in Aristotle's argument against the Platonic idealist conception, where he maintains that the unit is merely the measure of number. He also questions the equivalence of the point and the monas, but not the number :magnitude distinction. In this move away from a metaphysical, ontological foundation we see the beginning of the problematic of the status of number which became manifest in the Renaissance, eventually leading to a re-opening of the question as to the status of mathematical concepts generally.

When the question of the conceptualisation of number recurs in the Renaissance, we are at an extremely interesting historical juncture, since the problematic at this time could be seen as being both caused and resolved by pragmatism. As stated, fractions, as well as 'one' and zero had come to be regarded as numbers in everyday commercial usage, but this state of affairs was not so different from the everyday context of the Greek mathematical philosophers. The difference was that the Renaissance philosophies offered no perspective on questions of the conceptualisation of number, which only became problematic in connection with the algebraic developments.

In accordance with Diophantus' usage of the term 'arithmos' to represent the unknown in his investigations (where he also accepted fractional parts of the monas), it was assumed that the polynomial equations would yield numerical solutions. The strangeness of the emergent solutions, some of which were irreducible either to elements of the accepted number domain or to irrationals which corresponded to the Greek criterion of constructibility, caused a questioning of the earlier unselfconscious pragmatism and a concern to proscribe the limits of number. It is significant that the right of fractions to numberhood was not questioned.

In retrospect, we might say that the acceptance of fractions as numbers necessitated the eventual acceptance of the other candidates for numberhood: but that is to assume the concept of a number field, a number system which is closed under certain operations. This criterion, which is linked with a formalist understanding of mathematics, emerged very slowly. It was not until the 19th century that it came to be articulated more explicitly and receive more general acceptance, although it began to manifest in the consciously pragmatic attitudes adopted in the face of the new offspring of number. Until that time widely divergent views continued to be voiced as to the status of the negatives, imaginaries and irrationals; such views decreased in importance as the operational importance of the unplaced entities increased.

In the Renaissance climate of economic expansion (as opposed to the relatively static Greek economy) with its need for improved arithmetical techniques, and the belief in an effective, technologically oriented science based on the use of quantitative mathematical models, number was generally seen as quantity. The Greek ontological prohibition of fractions was no longer relevant, indeed it was scarcely considered;

fractions appeared 'natural' (the Greek view of them as ratios survived in the nomenclature 'rational' numbers); negatives were considered variously 'absurd', 'impossible', fictitious, and 'false', whilst the term 'imaginary' was coined dismissively by Descartes, complex numbers having been described earlier as 'useless' and 'sophistic'[49].

All these terms are revealing as regards the implied grounds for accepting candidates for numberhood; on the one hand there is an assumption of a sane, possible, real, true area prescribed by the accepted, positive, rational numbers (whose ontology is, at that time, not questioned); then, the term 'useless' betrays an incipient teleology in mathematics, it signals the beginnings of a more conscious pragmatism.

In the first place it is a contradictory pragmatism: Cardan denigrates imaginary numbers as 'useless' and continues to use them ('*putting aside the mental tortures involved*'[50]; theoretical attitudes oppose operational practice. The practice continues and in Girard's position we see the theoretical attitude reversed to form a coherent pragmatism: he asks himself:

Of what use are these impossible solutions [i.e. complex roots]? I answer: For three things - for the certitude of the general rules, for their utility, and because there are not other solutions.[51]

He frames his question pragmatically and answers himself partially tautologically; his further justification reveals a rudimentary formalism. From the perspective of modern rigour, this claim of certainty, and the assumption that the extant solutions complete the system, is unwarranted. The claim of certainty accompanied the use of the new numerical offspring throughout the long period of operation with them while their status was not agreed, not determined. It is only with the rebirth of rigour that the notion of certainty begins to become more specific. In the formalist doctrine, it finally returns to revive its original Greek counterpart: non-contradiction again becomes the paradigm, but the locus of fundamental validation is reduced; for the Greeks it was a consistent, global ontology; in the modern age completeness is required not in a global philosophy, rather it is sought in the local mathematical microcosm[52]: the original formalist demand being that a mathematical system be consistent and complete. Girard's theoretical attitude articulated the practice of the time; most mathematicians continued to practise pragmatism in contradiction with their theoretical positions[53].

The status of the irrationals was also re-evaluated. Some mathematicians (including Pascal and Barrow) wished to salvage the number:magnitude distinction although its original foundations were no longer valid ontologically or technically[54]. But there was also an attempt to articulate a new criterion for the acceptability of candidates for numberhood, more in keeping with the quantifying spirit of the Renaissance, based on the extension of decimal notation to include fractions: Stifel first gives the pragmatic reasons that might be given for accepting irrationals as numbers:

*Since, in proving geometric figures, when rational numbers fail us, irrational numbers take their place and prove exactly those things which rational numbers could not prove ... we are moved and compelled to assert that they truly are numbers, compelled that is, by the results which follow from their use - results*

---

49 Kline (1972), p.253
50 Kline (1972), p.253
51 Kline (1972), p.253
52 Once this reduction in philosophical demand is recognised, it can scarcely be viewed other than as an impoverishment.
53 Girard was an oddball in his time in his attempt to be consistent. Since then the split between doing and being has continued, deepened.
54 See Whiteside (1960-62)

*which we perceive to be real, certain and constant.*

But then he argues that:

*Other considerations compel us to deny that irrational numbers are numbers at all. To wit, when we seek to subject them to numeration* [i.e. decimal representation] *... we find that they flee away perpetually, so that not one of them can be apprehended precisely in itself ... Now that cannot be called a true number which is of such a nature that it lacks precision ... Therefore, just as an infinite number is not a number, so an irrational number is not a true number, but lies in a kind of a cloud of infinity*[55].

His subsequent argument that 'true' numbers are either whole numbers or fractions is circular, in that the decimal fractional notation is specifically constructed so as to express fractions in terms of ordered sequences of integers, and, in fact, contains a technical inconsistency in that he omits to point out the distinction between the ultimately predictable iterative procedure by which fractions whose denominators are prime to ten, when expressed in decimal notation, similarly recede to infinity, and the unpredictable course of the irrationals in such a representation. But in his attempt to abstract a criterion for acceptance of candidates to numberhood, that of precise locatability in terms of the extended decimal notation, he reveals an awareness of the problematic of the limit concept which is ignored in Stevin's later conception.

Ultimately Stifel's locatability criterion, like Descartes' argument for the limited acceptance of negative numbers (that equations with 'false', i.e. negative, roots can be transformed so as to yield positive roots), does not come to terms with the issue of the conceptualisation of number, but rather refers back to a <u>presupposed</u>, accepted number domain. In fact, Descartes does deal more directly with some important aspects of the number concept, but before discussing these we must first review Stevin's position.

Stevin's conception of his mathematical task was as part of a general project to recover the knowledge of a 'wise age' that he believed to have existed before the Greeks (whose culture he sees as the beginning of a 'barbarous age'). This compelled him to explicate thoroughly his radical conceptualisation of number, with reference to the surviving formulations of the Greeks.

Stevin was not the first to use the decimal fractional notation (Vieta and Bombelli amongst others used versions of such a notation) but his *Disme*, which functioned as a teaching manual for calculating with decimals and as propaganda advocating decimal standardisations of measures, was extremely influential[56] in spreading the use of the notation which clearly played a key role in determining his conception of number.

His first definition of number (articulating the spirit of the time) is *'that by which the quantity of each thing is expressed'*[57]. He then posits an analogy between number and continuous, homogeneous matter, stating his fundamental premise that the part is *'of the same material'* ('*de mesmematière*') as the whole[58]. Then, using the classical definition of number as a multitude of units, he argues that the unit, being a part of a multitude of units, is of the same material as the multitude of units, but the material of the multitude of units is 'number', so the material of the unit, and thus the unit itself, is 'number' - not the principle or beginning of number (as it is in the classical Greek view). He remarks by way of illustration that to deny this last step would be like denying that a piece of bread is bread; the devaluation of the integers (which he later makes explicit) is already clear. This conception of number does not allow for the qualitative

---

55 Kline (1972) p.251
56 The *Disme* was written in the vernacular. Stevin was one of the first to make this important change.
57 Klein (1968).
58 Klein (1968). It is telling that both Stevin and Descartes assume matter to be homogeneous.

difference between a loaf of bread and a pile of crumbs!

For Stevin, the unit is divisible (he invokes Diophantus); zero is now the beginning of number and the analogue of the geometric point. He states explicitly that *'number is by no means discontinuous'*. Number and magnitude are now so similar as to be almost identical; he attacks the use of the terms 'absurd' or 'irrational' for incommensurables: any root of a number is a number, since it is a part of a number. He makes a distinction between 'arithmetic' number and 'geometric' number, but this is actually a new distinction. 'Arithmetic' number is one expressed 'without an adjective of size' and 'geometric' numbers are 'quadratic', 'cubic' numbers etc. Any 'arithmetic' number may be a 'geometric' number, but when the numerical value is not known, 'geometric' numbers represent the indeterminate quantities in algebraic calculations, and are denoted by (1), (2), (3) etc. (where we would write powers of x, $x^2$, $x^3$ etc.). He thinks in termsof 'unknowns' which are still geometrically cloaked. The idea of variables enters mathematics later.

Stevin, with his practical background, was primarily interested in determinate solutions to problems; Vieta's approach was different. Like Stevin, he saw his work as a recovery of lost knowledge not as creation[59]. He saw in the study of polynomial equations the possibility of a general mathematical method. In his '*arsanalytice*'(the explicit telos of which was *'to leave no problem unsolved'*) he extended the Diophantine algebra, extracting general methods from Diophantus' particular cases (in this he saw himself as revealing the methods which Diophantus used but concealed). In the symbolic notation which he began to develop for polynomial equations the unknown is clearly distinguished by letter. He retains the terms 'side', 'square', 'cube', 'square-squared' verbally and refers to the quantities with which he deals, as 'magnitudes'. He does not use a sign for equality in the '*arsanalytice*'; he verbalises the progressive steps; it is not our modern equation form.

We are now in a position to return to Descartes (in whom we again find the belief in an earlier, more complete knowledge[60] ).Although the extent to which Descartes was influenced by Stevin and Vieta is not clear, his work could be seen, both on a mathematical, operational level, and on a more general philosophical level, as combining and extending theirs.

As regards mathematics 'proper', it is Descartes who liberates algebra from its internal geometricisation (which had become increasingly a relic from the Greek formulation) on two levels. Firstly, on an elemental, notational level, he synthesised Vieta's literal notation for the indeterminate term and Stevin's use of numbers to denote the power, thus creating the basis of our modern algebraic notation. Secondly, on an operational, conceptual level, he considered quantity as distinct from the geometric status (as Stevin had done for determinate number, but not in the case of the unknown): Descartes explicitly stated that a product of lines can be a line[61]. Then, having purged algebra of its geometric residue, he was able to establish a new relationship, at a higher structural level, between geometry and the algebra of dimensionless measure, a correspondence between equations and geometric curves; he formed a new synthesis, a new mode, coordinate geometry.

Descartes himself saw this as an example of the practical effectiveness of his general method which was aimed at the 'mathesisuniversalis', a general science of order and measurement, that could be seen as a further stage in achieving the generality which Vieta envisaged in his '*arsanalytice*'. Descartes founds his mathesis on a substantiation of Stevin's number-matter analogy; it is no longer seen as a metaphor; on the basis of his psycho-physiological model Descartes argues that there is an exact, real correspondence between number and matter; extension is both symbolic, as the object of general algebra, and real, as the substance of the corporeal world[62].

---

[59] Klein (1968).
[60] Klein (1968).
[61] Descartes (1968), p.169.

Descartes thus supplies both a philosophical and technical foundation for the budding science. His philosophy articulates the rationalist method and identifies the basic subject matter of abstract mathematics with that of science, the investigation of the material world. Technically his structural quantification of Euclidean space prepared the ground for the differential calculus, a mathematical technique for dealing with mechanical, quantitative change. A space ordered by a continuous (potentially infinitely divisible) measure was a necessary prerequisite for this development, as was a homogeneous time, a vitally important factor in the growth of Renaissance science. The homogeneity imposed upon space was imposed upon time[63]. Without the denial of the essential difference between these two basic orders of the life-world, the mechanical mathematics of the differential calculus is unthinkable. This is the crux of the Eleatic paradoxes. The Greek mathematical solution was to homogenise space and ignore time, deny change. The mathematical mode engendered could entertain an integral calculus, the method of exhaustion for accumulative approximation to the space contained by a curve (line or surface), but not investigation of a point phenomenon.

In the Renaissance, time enters as a conscious, explicit concern of science: the static, Greek 'episteme' (ἐπιστημη) gives way to Renaissance, time-and effect-oriented science. In mathematics Vieta's antinomy determinate:indeterminate replaces the discrete: continuous polarity which dominated Greek mathematics. The Greek conception of number was architecturally spatial; number was composed of geometrical arrangements of monads, indivisible units. The elements of arithmetike embodied spatial forms. Time is introduced into the conceptual material of mathematics (and number) through polynomial equations (as distinct from its explicit, external entry as an object of study). Originally they were seen as determinate; their subject was the 'unknown' which was specified geometrically. The change in perspective by which the 'unknown' became the 'variable' simultaneously recognised time and stepped outside of it; the determinate solution became subordinate to the vision of the form of the possibilities of solution.

This change in perspective is already inherent in the new notion of number as articulated by Stevin: number is homogeneous material, a conceptual object and the material which constitutes the object. The discrete:continuous contradistinction appears to be effaced. Stevin in fact argues the relativity of incommensurability on the grounds that a length that is incommensurable in one system, could be commensurable in another with a different standard unit; this fact, which is due to the relative nature of dimensional units, does not eradicate the contradiction inherent within number.

Stevin's formulation articulates the intuitive notion of his time; his concept of number is closer to the Greek magnitude, continuous measure (viz. Vieta's use of the term 'magnitude' for the subject of a polynomial equation); the absoluteness of the unit, the Greek foundation for discrete, heterogeneous number is undermined. The original integral decimal notation emphasised the homogeneous aspect of the integers, the repetitive procedure for approaching the potential infinity of succession. The extension of the notation in the other direction, i.e. to the infinitely small, gives an intuitive sense to the infinity of density, again emphasising the homogeneity. The decimal fractional notation, with its indefinitely close approximations to incommensurables encourages the illusion of a smooth elision from discrete to continuous, with the concomitant devaluation of the integers.

Descartes' explication of the latent idea of number-line completes the image. Decimals are originally constructed from integers; once they are constructed, the integers appear to be made up of decimals as (evanescent) building blocks or (a more coherent view permitted by Descartes' model) to be arbitrary places on a number line. In fact, at each decimal place one meets only another level of units; the place of transition from discrete to continuous continually recedes; one never confronts the essential difference between exact numbers and incommensurables; it is merely postponed indefinitely. Smooth number glosses over the integer:magnitude distinction; the Greek hypostatisation which had contained the discrete:continuous duality by keeping the two poles rigidly apart, is now dissolved. Its continuing relevance to considerations of the integers is obscured. It appears temporarily to have been banished; but it is merely transferred, transformed into an internal contradiction of the expanded number concept, which makes it possible for the calculus to approach the limit point. Under probing the contradiction again explodes.

---

[62] Klein (1968), pp.210-11.
[63] See Meyerson (1930) for an introduction to the deep cultural import of this scientific fiction ('fiction' in the technical sense as explicated in Vaihinger (1935)).

With the suppression of the discrete, the polarity undergoes a modal shift: the place formerly occupied by the discrete, the number object, now houses determinate, finite number, any specific representative of number the material, the continuous, which now represents the potentially infinite divisibility and/or range possible for any determinate number; it is simultaneously the material from which any determinate number is formed and the homogeneous number line (or space) from which a determinate number may be chosen. Stevin states that every 'arithmetic' (determinate) number is the beginning of 'geometric' (indeterminate) number, just as zero is the beginning of 'arithmetic' number[64].

Whereas the Greek number :magnitude distinction was a static horizontal dualism fixed in space, the new concept of number implicitly contains the notion of a variable. The polarity operates between the levels, determinate and indeterminate; it is vertical rather than horizontal, dialectical rather than static. The new model fuses space and time in common homogeneity: determinate number is both formed of, and chosen from homogeneous, potential number; the Greek spatial composition of number persists in the identification of a determinate number with a line segment; but a determinate number is also a place, a point on a number line which may be singled out like a moment of time from which the probable past and the possible future stretch endlessly away.

The implications of this mode of fusion are immense and warrant further discussion, but immediately we see the root of the problematic that concerned Miller and Meyerson[65]; the new mathematics initiated in the Renaissance deals with time by reducing it to quantified space. The past and future are arbitrary in this schema (unlike the life-world), dependent on an arbitrary origin; the illusion is created of the possibility of stepping outside time, an illusion which, like the fiction of Euclidean space itself, is valid only within certain limits. A framework is created for the quantitative description (and thus possible prediction) of mechanical change, but the essence of experienced time, real change, the emergence of the qualitatively new, still eludes description.

This problematic has recently been approached within mathematics by Rene Thom[66]; his catastrophe theory can provide qualitative models for changes of state within a certain nexus, but he himself considers that the theory is essentially incapable of adaptation to quantification for predictive scientific purposes. The Renaissance vision of knowledge that is both certain and effective is reaching the limits of the mode it engendered.

It would be worthwhile to return to examine more closely the roots of present mathematical problematics in the Renaissance incunabular of our mathematical mode, as well as the continuities and reversals from the original Greek seeds. There are, in fact, several interesting parallels between the Greek beginnings and the Renaissance rebirth.

---

64 Klein (1968). I think that further investigation of this profound, conceptual shift could yield greater understanding not only of the ontology of numbers and mathematics but also of the underlying consciousness of our mathematico- scientific culture.
65 See Miller (1948) and Meyerson (1930).
66 See Thom (1975)

The Pythagorean mathematical inspiration was the vision of a numerical description of the world, where 'number' was discrete, heterogeneous and in some way material. This thesis met with contradiction in the form of the irrationals; the idea of number as discrete was not sacrificed (that would have meant abdication of their metaphysics as a changeless reality}; but a new rigour was necessary to ensure non-contradiction (equated with certainty). On the one hand it was necessary to supply a foundation for operation with irrationals; on the other, the general form of mathematical demonstration was tightened. In the Renaissance twist of the spiral, there is still a vision of a numerical world-description ('number' is now abstract and continuous, having subsumed the discrete:continuous antinomy) but with a new impetus deriving from the perception of mathematics as an epistemological method, a general art for solving problems.

The Pythagorean vision is of an isomorphism between the realm of discrete number and the life-world; the demonstrative method is secondary. Descartes reiterates this vision, with continuous number replacing the discrete, interposing a level of abstraction; he reverses priorities: his epistemological, rationalist method, modelled on the mathematical proof form, is primary. He also reverses the role of the form: in its original context the theorem precedes the proof, its direct function is static (it is only through a further act of reflection on the internal structural components of the proof, an involution, in accordance with Lakatos' proof analysis, that it receives a function in generating new understandings); when Descartes appropriates the form he interprets it as an epistemological method whose function is generative[67]. (This part of his dream was never realised) In these reversals, the monadic structure retains its primacy. For the Pythagorean, the monadology is ontological; the unit elements are the prime constituents of the phenomenal world. For Descartes, it is epistemological: the constituents are the clear and distinct ideas which accumulate to form a body of sure knowledge.

There is a further parallel, both with regard to content and history, between the (Pythagorean) Platonic metaphysics and the Cartesian rationalist science. Both deny real, qualitative change, the first by denying time, the second, by neutering it.. Both proved successful within limits; their very success caused the ideas to become embedded as ideology, so that even when the limits of their validity are approached and the original doctrines are questioned, their consequences still survive in mathematical praxis of which the roots have been forgotten, sedimented in history.

The rebirth of Greek mathematics in Western culture is simultaneously a completion and an inversion. The development of the decimal notation fosters the new concept of homogeneous number, allowing time to be subsumed into mathematics. The initial consistent place system, of course, depends on the existence of a symbol for zero, and it would be worthwhile considering the further implications of this, for instance, the change of meaning that occurred in the cultural transfer of zero from its Hindu origins, where there is a sense of a full nothingness, to a society where 'nothing' is a mere absence, where 'nature abhors a vacuum' etc. What are the implications of our conception of zero and those of Hindu and Buddhist philosophy? Such questions obviously relate to considerations of the calculus. Interestingly Buddhist logic, beginning from the point-instant as the basic reality, i.e. almost the opposite metaphysics to the Platonist, came very close to notions of a differential calculus[68].

The mode of the symbolic notation developed in the Renaissance could obviously be examined in more detail. For instance, before the Renaissance, the equation was not the paradigm of mathematics. The primacy of the equation (which has only very recently cometo be questioned) internalises the mode of the axiomatic method: static consolidation of atemporal positivities[69]. The question of overdetermination of mathematical symbols needs to be investigated.

---

67 As stated this seems to have been more of a misunderstanding on Descartes' part than a deliberate step.
68 See Stcherbatsky (1962)

When we look at mathematics as a language, we see that the concept of number which emerges in the Renaissance is adjectival with respect to ordinary language. The number language then bears a skew relation to ordinary language, since these adjectives are then accorded a substantive function in mathematical grammar. The situation is, of course, vastly more complex; mathematical grammar is not isomorphic to the grammar of ordinary language, but the attempt to understand mathematics in this way is valuable. The divergence between mathematical language and ordinary language only began with the development of symbolic notation in the Renaissance; before that time mathematics was still predominantly verbal.

In classical Greece mathematical language was embedded in ordinary language: to understand mathematical objects was to locate them in a global ontology; numbers attained reality as objects, and so could coherently be understood as substantives, by embodying geometric forms. In the Renaissance and Enlightenment mathematical operations are increasingly symbolised. Mathematical language is thus formally separated from verbal language. It has an internal coherence, being now composed homogeneously of symbolic elements. It is free to follow its own dynamic according to its own, internal, grammatical and syntactical laws. In the new reflexivity of mathematics the number concept is extended by internalising mathematical operations; an operational classification of number supersedes the Greek geometric classification.

In the initial phases of growth of this mode as Klein says:

*the whole complex of ontological problems which surrounds the ancient concept of number loses its object in the context of the symbolic conception, since there is no immediate occasion for questioning the mode of being of the 'symbol' itself.[70]*

The extent of the rupture was not, at first, recognised; Descartes' identification of the mathematical object-world with the perceived material world and Kant's attempt to refound it metaphysically, was not called into question until the development of consistent non-Euclidean geometries, causing attempts to justify what had by that time become the status quo: conventionalism, logicism, formalism. These explicitly refuse to consider the problematic of the relation between mathematics and the life-world, attempting instead to create a Frankenstein's monster of mathematics, a self-sufficient entity whose judgment is more certain than that of its creators. The different ways in which the various attempts failed to achieve their original stated aims, merits further examination, but the very fact of the failures calls for a re-examination of the basis for such attempts. When Wittgenstein questions the motivation behind the foundations fervour, he asks:

*But what was the attempt made for? Was it not due to an uncertainty in another place?[71]*

---

69 See Adorno and Horkheimer (1944), p.7 ff, for a description of the mutual mirroring of mathematical and societal developments, for example, "Bourgeois society is ruled by equivalence. It makes the similar comparable by reducing it to abstract quantities. To the Enlightenment, that which does not reduce to numbers, and ultimately to the one, becomes illusion; modern positivism writes it off as literature. Unity is the slogan from Parmenides to Russell. The destruction of gods and qualities alike is insisted upon."
70 Klein (1968).
71 Wittgenstein (1956)

Is that 'other place' not the questionable place of mathematics itself in the life-world, entailing consideration of such problematic concepts as certainty which bridge the objective and subjective? The problematic posed for the Greeks by the symbolic status of number, as to its ontological reality, was temporarily submerged as the whole of mathematics assumed a symbolic character. Thus, the problematic now recurs on a larger scale, of the ontological status of mathematics.

## 5. SUMMARY

We have seen that at its inception in classical Greek times, mathematics was a holistic, practical philosophy, concerned not only with technical, quantitative knowledge of the physical world but also with a qualitative understanding of the nature of knowledge and human life as a whole. Arithmetike, the study of the qualities of number (i.e. arithmoi, the natural numbers) was central to this.

The discovery of the existence of irrational numbers posed problems for this global mathematics, which were not met explicitly. Over the following centuries, the distinction between number (the discrete, positive integers) and magnitude (continuously divisible, physical size) was increasingly ignored with the technical expansion of mathematics. Mathematical pragmatism began to supplant the Greek mathematical ontology.

The 13th century introduction into Europe through commerce of the Hindu-Arabic decimal numeral system, replacing the cumbersome Roman numerals, made possible the emergence of the revolutionary idea of a number line. Together with rapidly developing mathematical symbolism, this culminated, in the 17th century, with Cartesian algebraic geometry and the differential calculus, which have dominated mathematical activity to a large extent since. We saw that the changes here in the concepts of number and mathematics generally were intimately connected with changes in the concept of time. Equations were no longer assumed to be determinate: the notion of the variable was born, The number line as homogeneous, infinitely divisible measure, not only ordered space; it was also imposed upon time. This allowed the extremely powerful, quantitative description of mechanical change but denied time's essentially different nature. The static Platonic (Pythagorean) metaphysics was replaced by causal science ushering in Bacon's innovative notion of technological progress.

We began to disentangle some of the far-reaching implications involved in the intricate interrelations of the mathematical conceptual changes in this period. Mathematical symbolism reached such a degree of complexity that its vocabulary, grammar and syntax parted company with its source verbal matrix. Much work remains to be done here since this abstract mathematics where complex concepts are locked inside seemingly simple signs, such as those for zero, infinity, 'equals', variables, functions etc, is a deep part of our cultural, sedimented history.

By the end of the 19th century contradictions again arose. Formalism was an attempt to establish certainty in mathematics on the basis of the logical proof form. The attempt failed in the face of Gödel's theorems. Wittgenstein pointed out that the desire for foundations can not be satisfied within mathematics. Its philosophy needs to be based in the larger context of the life-world. We investigated possibilities for such a philosophy in phenomenology and Lakatos' methodology of mathematical development, given that the latter is inherent in mathematical history, waiting to be consciously adopted. Such approaches recognise mathematics as a living, creative process, not just product.

---

This could be the next stage of mathematics' relationship with time: neither denying it as in Platonic ontology, nor neutering it as in Cartesian rationalist science, but recognising that mathematics itself is situated within historical time. Becoming self-reflective it could claiming its history, and make choices on the basis of this recognition. One choice might be to investigate the possibilities of recognising new integer qualities in the wonderful mathematical world of the 21$^{st}$ century: the specificities of different n-dimensional spaces, for example.

One might say that with the loss of importance of the integers, mathematics lost its integrity. Perhaps restoring them to their rightful status (in accordance with Gauss' view of number theory as the queen of mathematics) might be a turn of the spiral whereby mathematics discovers or recreates a new integrity.

**ACKNOWLEDGEMENTS**

I thank my father, Bill Graves, for his unquestioning belief in and support for my work and my mother, Brenda Graves, for her desire that I fulfil myself. I thank my husband Jon Gregory, for his tolerance, my daughter, Jemima Gregory, for her enlightening being and the child in my womb for supplying me with an irrevocable deadline, which is naturally a birth-time. I thank Helle Munro for her commitment to my boundary conditions.

In the academic field, I thank the following for their encouragement and/or criticism: (in chronological order) Eleanor Gill, John Conway, Ralph Schwarzenberger, David Fowler, SigurdZienau, Nicholas Maxwell, P.M. Rattansi, Barry Smith, Kevin Mulligan, Clive Kilmister, Caroline Dunmore, Peter Saunders and Alison Watson.

AmoAmitabhaya.